\documentclass[english,oneside]{amsart}
\usepackage[T1]{fontenc}
\usepackage[latin1]{inputenc}
\usepackage{amssymb}

\makeatletter
 \theoremstyle{plain}
\newtheorem{thm}{Theorem}[section]
  \theoremstyle{plain}
  \newtheorem{lem}[thm]{Lemma}
  \theoremstyle{definition}
  \newtheorem{defn}[thm]{Definition}
  \theoremstyle{remark}
  \newtheorem{rem}[thm]{Remark}

\usepackage{babel}
\makeatother
\begin{document}

\title[Asymptotics for Stieltjes-Wigert Polynomials]{Plancherel-Rotach Asymptotics for Stieltjes-Wigert Orthogonal Polynomials
with Complex Scaling}

\curraddr{School of Mathematics\\
 Guangxi Normal University\\
 Guilin City, Guangxi 541004\\
 P. R. China.}

\author{Ruiming Zhang}

\email{ruimingzhang@yahoo.com}



\begin{abstract}
In this work we study the Plancherel-Rotach type asymptotics for Stieltjes-Wigert
orthogonal polynomials with complex scaling. The main term of the
asymptotics contains Ramanujan function $A_{q}(z)$ for the scaling
parameter on the vertical line $\Re(s)=2$, while the main term of
the asymptotics involves the theta functions for the scaling parameter
in the strip $0<\Re(s)<2$. In the latter case the number theoretical
property of the scaling parameter completely determines the order
of the error term. $ $ These asymptotic formulas may provide some
insights to some new random matrix model and also add a new link between
special functions and number theory.
\end{abstract}
\maketitle

\section{Introduction\label{sec:Introduction}}

The Plancherel-Rotach type asymptotics for classical orthogonal polynomials
are essential to obtain universality results in random matrix theory.
The associated random matrix models for $q$-orthogonal polynomials
are still unknown today. It might be interesting to calculate the
Plancherel-Rotach type asymptotics for $q$-orthogonal polynomials
to gain some insights to the related random matrix models.

In \cite{Ismail7} we studied Plancherel\emph{-}Rotach type asymptotics
for three families of $q$-orthogonal polynomials with real logarithmic
scalings. They are Ismail-Masson polynomials $\left\{ h_{n}(x|q)\right\} _{n=0}^{\infty}$
, Stieltjes-Wigert polynomials $\left\{ S_{n}(x;q)\right\} _{n=0}^{\infty}$
and $q$-Laguerre polynomials $\left\{ L_{n}^{(\alpha)}(x;q)\right\} _{n=0}^{\infty}$.
They are $q$-orthogonal polynomials associated with indeterminant
moment problems. The asymptotics reveal a remarkable pattern which
is quite different to the pattern associated with the classical Plancherel\emph{-}Rotach
asymptotics\cite{Szego,Ismail2}. The main term of asymptotics may
contain Ramanujan function $A_{q}(z)$ or theta function according
to the value of the scaling parameter.

In this paper, we study the Plancherel-Rotach type asymptotics for
Stieltjes-Wigert polynomials under complex scaling. We also have explicit
error bounds for the asymptotic formulas. When the scaling parameter
$s$ is in the vertical strip $0<\Re(s)<2$, the order of the error
term is completely determined by the number theoretical property of
the scaling parameter $s$.

In \S\ref{sec:Preliminaries}, we list some facts from $q$-series
and number theory. We present the asymptotic formulas for Stieltjes-Wigert
polynomials under complex scaling in \S\ref{sec:Main Results}. Finally
we use a discrete analogue of Laplace method to derive these formulas
in \S\ref{sec:Stieltjes-Wigert}.

Through out this paper, We shall always assume that $0<q<1$ unless
otherwise stated. We also assume that $s=\sigma+it$ is a complex
number with $\sigma=\Re(s)$ and $t=\Im(s)$. All the $\log$ and
power functions are taken as their respective principle branches.

\section{Preliminaries\label{sec:Preliminaries}}

\subsection{Something from $q$-series}

In this paper we follow the usual notation for $q$-series, \cite{Andrews4,Gasper,Ismail2}

\begin{equation}
(a;q)_{0}:=1\quad(a;q)_{n}:=\prod_{k=0}^{n}(1-aq^{k}),\quad\left[\begin{array}{c}
n\\
k\end{array}\right]_{q}:=\frac{(q;q)_{n}}{(q;q)_{k}(q;q)_{n-k}},\label{eq:2.1}\end{equation}
 and $n=\infty$ is allowed in the above definitions. Assume that
$|z|<1$, the $q$-Binomial theorem is \cite{Andrews4,Gasper,Ismail2}
\begin{equation}
\frac{(az;q)_{\infty}}{(z;q)_{\infty}}=\sum_{k=0}^{\infty}\frac{(a;q)_{k}}{(q;q)_{k}}z^{k},\label{eq:2.2}\end{equation}
 which defines an analytic function in the region $|z|<1$. Its limiting
case, \begin{equation}
(z;q)_{\infty}=\sum_{k=0}^{\infty}\frac{q^{k(k-1)/2}}{(q;q)_{k}}(-z)^{k}\quad z\in\mathbb{C}\label{eq:2.3}\end{equation}
 is one of many $q$-exponential identities. Ramanujan function $A_{q}(z)$
is defined as \cite{Ramanujan,Ismail2} \begin{equation}
A_{q}(z):=\sum_{k=0}^{\infty}\frac{q^{k^{2}}}{(q;q)_{k}}(-z)^{k}.\label{eq:2.4}\end{equation}
 In order to express our error terms concisely, we also define a related
function $B_{q}(z)$, \begin{equation}
B_{q}(z):=\sum_{k=0}^{\infty}\frac{q^{k^{2}}z^{k}}{(q;q)_{k}}.\label{eq:2.5}\end{equation}
 Clearly,\begin{equation}
\left|A_{q}(z)\right|\le B_{q}(|z|).\label{eq:2.6}\end{equation}
 Since\[
k+1\ge\frac{1-q^{k+1}}{1-q}\ge(k+1)q^{k}\]
 and\begin{eqnarray*}
{B'}_{q}(|z|) & = & q\sum_{k=0}^{\infty}\frac{q^{k^{2}}|qz|^{k}(k+1)q^{k}}{(q;q)_{k}(1-q^{k+1})}\\
 & \le & \frac{q}{1-q}B_{q}(|z|),\end{eqnarray*}
 thus,\begin{equation}
\left|A'_{q}(z)\right|\le{B'}_{q}(|z|)\le\frac{q}{1-q}B_{q}(|z|).\label{eq:2.7}\end{equation}
 The theta function, \cite{Andrews4,Gasper,Ismail2}\begin{equation}
\Theta(z|q):=\sum_{n=-\infty}^{\infty}q^{n^{2}}z^{n}\label{eq:2.8}\end{equation}
 is clearly defined for any complex number $z\neq0$. Jacobi triple
product identity is\begin{equation}
\Theta(z|q)=(q^{2},-qz,-q/z;q^{2})_{\infty}.\label{eq:2.9}\end{equation}
 Stieltjes\emph{-}Wigert orthogonal polynomials $\left\{ S_{n}(x;q)\right\} _{n=0}^{\infty}$
are defined as \cite{Ismail2}\begin{equation}
S_{n}(x;q)=\sum_{k=0}^{n}\frac{q^{k^{2}}(-x)^{k}}{(q;q)_{k}(q;q)_{n-k}}.\label{eq:2.10}\end{equation}
 Let \begin{equation}
s=\tau+2+i\frac{2\theta\pi}{\log q}\quad\tau,\theta\in\mathbb{R},\label{eq:2.11}\end{equation}
 and\begin{equation}
x_{n}(z,s)=zq^{-ns}.\label{eq:2.12}\end{equation}
 Reverse the summation order of \eqref{eq:2.10} to get\begin{equation}
\frac{S_{n}(x_{n}(z,s);q)}{(-z)^{n}q^{n^{2}(1-s)}}=\sum_{k=0}^{n}\frac{q^{k^{2}}e^{2nk\theta\pi i}}{(q;q)_{k}(q;q)_{n-k}}\left(-\frac{q^{\tau n}}{z}\right)^{k}.\label{eq:2.13}\end{equation}
 For any real number $x$, then, \begin{equation}
x=\left\lfloor x\right\rfloor +\left\{ x\right\} ,\label{eq:2.14}\end{equation}
 where the fractional part of $x$ is $\left\{ x\right\} \in[0,1)$
and $\left\lfloor x\right\rfloor \in\mathbb{Z}$ is the greatest integer
less or equal $x$. The arithmetic function \begin{equation}
\chi(n)=2\left\{ \frac{n}{2}\right\} =n-2\left\lfloor \frac{n}{2}\right\rfloor ,\label{eq:2.15}\end{equation}
 which is the principal character modulo $2$,\begin{equation}
\chi(n)=\begin{cases}
1 & 2\nmid n\\
0 & 2\mid n\end{cases}.\label{eq:2.16}\end{equation}
 We will also make uses of the trivial inequality \begin{equation}
|e^{z}-1|\le|z|e^{|z|}\label{eq:2.17}\end{equation}
 for any $z\in\mathbb{C}$. The following lemma is from \cite{Ismail7}.

\begin{lem}
\label{lem:1}Given any $n\in\mathbb{N}$, if $a>0$, \begin{equation}
\frac{(a;q)_{\infty}}{(a;q)_{n}}=(aq^{n};q)_{\infty}=1+R_{1}(a;n)\label{eq:2.18}\end{equation}
 with\begin{equation}
\left|R_{1}(a;n)\right|\le\frac{(-aq^{2};q)_{\infty}aq^{n}}{1-q},\label{eq:2.19}\end{equation}
 while for $0<aq<1$, \begin{equation}
\frac{(a;q)_{n}}{(a;q)_{\infty}}=\frac{1}{(aq^{n};q)_{\infty}}=1+R_{2}(a;n)\label{eq:2.20}\end{equation}
 with \begin{equation}
\left|R_{2}(a;n)\right|\le\frac{aq^{n}}{(1-q)(aq;q)_{\infty}}.\label{eq:2.21}\end{equation}

\end{lem}
\begin{proof}
It is clear from \eqref{eq:2.3} that\[
R_{1}(a;n)=\sum_{k=1}^{\infty}\frac{q^{k(k-1)/2}(-aq^{n})^{k}}{(q;q)_{k}}=-aq^{n}\sum_{k=0}^{\infty}\frac{q^{k(k+1)/2}(-aq^{n})^{k}}{(q;q)_{k+1}}.\]
 Hence for $a>0$, we have \begin{eqnarray*}
\left|\sum_{k=0}^{\infty}\frac{q^{k(k+1)/2}(-aq^{n})^{k}}{(q;q)_{k+1}}\right| & \le & \frac{1}{1-q}\sum_{k=0}^{\infty}\frac{q^{k(k-1)/2}}{(q;q)_{k}}(aq^{n+1})^{k}\\
 & \le & \frac{(-aq^{2};q)_{\infty}}{1-q},\end{eqnarray*}
 and \eqref{eq:2.19} follows. Moreover from \eqref{eq:2.2} \begin{eqnarray*}
R_{2}(a;n) & = & aq^{n}\sum_{k=0}^{\infty}\frac{\left(aq^{n}\right)^{k}}{(q;q)_{k+1}},\end{eqnarray*}
 hence \begin{eqnarray*}
\left|R_{2}(a;n)\right| & \le & \frac{aq^{n}}{(1-q)(aq^{n};q)_{\infty}}\le\frac{aq^{n}}{(1-q)(aq;q)_{\infty}},\end{eqnarray*}
 which is \eqref{eq:2.21} and the proof of the lemma is complete. 
\end{proof}
\[
\]

\subsection{Generalized Irrational Measure}

For an irrational number $\theta$, Chebyshev's theorem implies that
for any real number $\beta$ , there exist infinitely many pairs of
integers $n$ and $m$ with $n>0$ such that\cite{Hua} \begin{equation}
n\theta=m+\beta+a_{n}\quad{\rm with}\quad|a_{n}|\le\frac{3}{n}.\label{eq:2.22}\end{equation}
 Clearly, Chebyshev's theorem says that the arithmetic progression
$\left\{ n\theta\right\} _{n\in\mathbb{Z}}$ is ergodic in $\mathbb{R}$.

\begin{defn}
\label{def:generalized irrationality measure}Given real numbers $\theta_{j}$
and $\beta_{j}$ for $j=1,...,N$, the generalized irrationality measure
$\omega(\theta_{1},\dots,\theta_{N}|\beta_{1},\dots,\beta_{N})$ of
$\theta_{1},\dots,\theta_{N}$ associated with $\beta_{1},\dots,\beta_{N}$
is defined as the least upper bound of the set of real numbers $r$
such that there exist infinitely many integers $n$ and $m_{1},\dots,m_{N}$
with $n>0$ such that\begin{equation}
|n\theta_{j}-\beta_{j}-m_{j}|<\frac{1}{n^{r-1}}\label{eq:2.23}\end{equation}
 for $j=1,\dots,N$. 
\end{defn}
\begin{thm}
Let $\theta$ be an irrational number, then for any real number $\beta$,
its generalized irrational measure associated with $\beta$ is $\omega(\theta|\beta)\ge2$
. 
\end{thm}
\begin{proof}
This is a direct consequence of the Chebyshev's theorem. 
\end{proof}
The irrationality measure (or Liouville\emph{-}Roth constant) $\mu(\theta)$
of a real number $\theta$ is defined as the least upper bound of
the set of real numbers $r$ such that\cite{Wikipdedia} \begin{equation}
0<|n\theta-m|<\frac{1}{n^{r-1}}\label{eq:2.24}\end{equation}
 is satisfied by an infinite number of integer pairs $(n,m)$ with
$n>0$. It is clear that we have \begin{equation}
\omega(\theta|0)=\mu(\theta).\label{eq:2.25}\end{equation}
 A real algebraic number $\theta$ of degree $l$ if it is a root
of an irreducible polynomial of degree $l$ with integer coefficients.
Liouville's theorem in number theory says that for a real algebraic
number $\theta$ of degree $l$, there exists a positive constant
$K(\theta)$ such that for any integer $m$ and $n>0$ we have \begin{equation}
|n\theta-m|>\frac{K(\theta)}{n^{l-1}}.\label{eq:2.26}\end{equation}
 A Liouville number is a real number $\theta$ such that for any positive
integer $l$ there exist infinitely many integers $n$ and $m$ with
$n>1$ such that\cite{Wikipdedia}\begin{equation}
0<|n\theta-m|<\frac{1}{n^{l-1}},\label{eq:2.27}\end{equation}
 and the terms in the continued fraction expansion of every Liouville
number are unbounded. Even though the set of all Liouville numbers
is of Lebesgue measure zero, it is known that almost all real numbers
are Liouville numbers topologically.

\begin{thm}
The generalized irrational measure $\omega(\theta|\beta)$ has the
following properties: 
\begin{enumerate}
\item For any real algebraic number $\theta$ of degree $l$, one has $\omega(\theta|0)\le l$; 
\item For a Liouville number $\theta$, one has $\omega(\theta|0)=\infty$. 
\end{enumerate}
\end{thm}
\begin{proof}
The first assertion follows from the definition\ref{def:generalized irrationality measure}
and \eqref{eq:2.26} while the last assertion follows directly from
the definitions of the generalized irrational measure and the Liouville
numbers. 
\end{proof}
It is clear that the quadratic irrationals $\theta$ such as $\sqrt{2}$
have $\omega(\theta|0)=2$.

\section{Main Results\label{sec:Main Results} }

Given a real number $\theta$, we define\begin{equation}
\mathbb{S}(\theta)=\left\{ \left\{ n\theta\right\} :n\in\mathbb{N}\right\} .\label{eq:3.1}\end{equation}
 When $\theta$ is a rational number, $\mathbb{S}(\theta)$ is a finite
subset of $[0,1)$. When $\theta$ is irrational, the set $\mathbb{S}(\theta)$
is a dense subset of $(0,1)$, which is a consequence of Chebyshev's
theorem.

\begin{thm}
\label{thm:stieltjes-wigert}Given any nonzero complex number $z$,
let $s$ and $x_{n}(z,s)$ be defined as in \eqref{eq:2.15} and \eqref{eq:2.12},
we have the following results for Stieltjes-Wigert polynomials: 
\begin{enumerate}
\item When $\tau>0$, we have \begin{equation}
\frac{S_{n}(x_{n}(z,s);q)(q;q)_{\infty}}{(-z)^{n}q^{n^{2}(1-s)}}=1+r_{sw}(n|1)\label{eq:3.2}\end{equation}
 with\begin{equation}
\left|r_{sw}(n|1)\right|\le\frac{qB_{q}(q^{2}|z|^{-1})}{(1-q)|z|}.\label{eq:3.3}\end{equation}

\item When $\tau=0$ and $\theta$ is a rational number. For any $\lambda\in\mathbb{S}(\theta)$,
there are infinitely many pair of integers $n$ and $m$ with $n>0$
such that \begin{equation}
n\theta=m+\lambda,\quad\lambda=\left\{ n\theta\right\} .\label{eq:3.4}\end{equation}
 For each of such integers, \begin{equation}
\frac{S_{n}(x_{n}(z,s);q)(q;q)_{\infty}}{(-z)^{n}q^{n^{2}(1-s)}}=A_{q}\left(\frac{e^{2\lambda\pi i}}{z}\right)+r_{sw}(n|2)\label{eq:3.5}\end{equation}
 and the error is majorized as \begin{equation}
|r_{sw}(n|2)|\le\frac{2(-q^{3};q)_{\infty}B_{q}(|z|^{-1})}{(q;q)_{\infty}}\left\{ q^{n/2}+\frac{q^{n^{2}/4}}{|z|^{\left\lfloor n/2\right\rfloor }}\right\} .\label{eq:3.6}\end{equation}
 for $n$ is sufficiently large. 
\item When $\tau=0$ and $\theta$ is an irrational number. Given any real
number $\beta\in[0,1)$ and a real number $\rho$ with \begin{equation}
0<\rho<\omega(\theta|\beta)-1,\label{eq:3.7}\end{equation}
 there are infinitely many pair of integers $n$ and $m$ with $n>0$
such that\begin{equation}
n\theta=m+\beta+\gamma_{n}\quad|\gamma_{n}|\le\frac{1}{n^{\rho}}.\label{eq:3.8}\end{equation}
 Then \begin{equation}
\frac{S_{n}(x_{n}(z,s);q)}{(-z)^{n}q^{n^{2}(1-s)}}=A_{q}\left(\frac{e^{2\pi i\beta}}{z}\right)+e_{sw}(n|3)\label{eq:3.9}\end{equation}
 with\begin{equation}
|e_{sw}(n|3)|\le\frac{24(-q^{3};q)_{\infty}B_{q}(|z|^{-1})}{(q;q)_{\infty}}\left\{ \frac{\log^{2}n}{n^{\rho}}+\frac{q^{\nu_{n}^{2}}}{|z|^{\nu_{n}}}\right\} \label{eq:3.10}\end{equation}
 and\begin{equation}
\nu_{n}=\left\lfloor \frac{q^{4}\log^{2}n}{1+\log q^{-1}}\right\rfloor \label{eq:3.11}\end{equation}
 for $n$ is sufficiently large. 
\item When $-2<\tau<0$ and both $\tau$ and $\theta$ are rational. Assume
that for some $\lambda\in\mathbb{S}(-\tau)$ and some $\lambda_{1}\in\mathbb{S}(\theta)$,
there are infinite number of positive integers $n$ such that\begin{equation}
-n\tau=m+\lambda,\quad m=\left\lfloor -\tau n\right\rfloor ,\quad\lambda=\left\{ -\tau n\right\} \label{eq:3.12}\end{equation}
 and\begin{equation}
n\theta=m_{1}+\lambda_{1},\quad\lambda_{1}=\left\{ n\theta\right\} .\label{eq:3.13}\end{equation}
 Then for each such $n$, $m_{1}$, and $m_{2}$ , we have\begin{eqnarray}
S_{n}(x_{n}(z,s);q) & = & \frac{(-z)^{n}q^{n^{2}(1-s)+\left\lfloor m/2\right\rfloor \left(\tau n+\left\lfloor m/2\right\rfloor \right)}}{(q;q)_{\infty}^{2}(-ze^{-2n\theta\pi i})^{\left\lfloor m/2\right\rfloor }}\nonumber \\
 & \times & \left\{ \Theta\left(-zq^{\chi(m)+\lambda}e^{-2\pi i\lambda_{1}}\mid q\right)+r_{sw}(n|4)\right\} \label{eq:3.14}\end{eqnarray}
 with\begin{equation}
|r_{sw}(n|4)|\le\frac{6(-q^{3};q)_{\infty}^{2}\Theta(|z|\mid\sqrt{q})}{(1-q)^{2}}\left\{ q^{\nu_{n}}+q^{\nu_{n}^{2}}|z|^{\nu_{n}}+\frac{q^{\nu_{n}^{2}/2}}{{|z|}^{\nu_{n}}}\right\} \label{eq:3.15}\end{equation}
 and\begin{equation}
\nu_{n}=\min\left\{ \left\lfloor \frac{(2+\tau)n}{8}\right\rfloor ,\left\lfloor \frac{-\tau n}{8}\right\rfloor \right\} .\label{eq:3.16}\end{equation}
 for $n$ sufficiently large. 
\item When $-2<\tau<0$ , $\tau$ is rational and $\theta$ is irrational.
Given a real number $\beta\in[0,1)$ let \begin{equation}
0<\rho<\omega(\theta|\beta)-1,\label{eq:3.17}\end{equation}
 there are infinitely many positive integers $n$, \begin{equation}
n\theta=m_{1}+\beta+b_{n},\quad{|b}_{n}|\le\frac{1}{n^{\rho}}\label{eq:3.18}\end{equation}
 and some $\lambda\in\mathbb{S}(-\tau)$ with\begin{equation}
-n\tau=m+\lambda,\quad m=\left\lfloor -n\tau\right\rfloor .\label{eq:3.19}\end{equation}
 For each such $\lambda,\beta,n,m$, we have\begin{eqnarray}
S_{n}(x_{n}(z,s);q) & = & \frac{(-z)^{n}q^{n^{2}(1-s)+\left\lfloor m/2\right\rfloor \left(\tau n+\left\lfloor m/2\right\rfloor \right)}}{(q;q)_{\infty}^{2}(-ze^{-2n\theta\pi i})^{\left\lfloor m/2\right\rfloor }}\nonumber \\
 & \times & \left\{ \Theta\left(-zq^{\chi(m)+\lambda}e^{-2\beta\pi i}\mid q\right)+r_{sw}(n|5)\right\} \label{eq:3.20}\end{eqnarray}
 with\begin{eqnarray}
|r_{sw}(n|5)| & \le & \frac{54(-q^{3};q)_{\infty}^{2}\Theta(|z|\mid\sqrt{q})}{(1-q)^{2}}\nonumber \\
 & \times & \left\{ |z|^{\nu_{n}}q^{\nu_{n}^{2}}+\frac{q^{\nu_{n}^{2}/2}}{|z|^{\nu_{n}}}+\frac{\log^{2}n}{n^{\rho}}\right\} \label{eq:3.21}\end{eqnarray}
 and\begin{equation}
\nu_{n}=\left\lfloor \frac{q^{4}\log^{2}n}{1+\log q^{-1}}\right\rfloor \label{eq:3.22}\end{equation}
 for $n$ sufficiently large. 
\item When $-2<\tau<0$, $\tau$ is irrational and $\theta$ is rational,
Given a real number $\beta\in[0,1)$ let \begin{equation}
0<\rho<\omega(\theta|\beta)-1,\label{eq:3.23}\end{equation}
 there are infinitely many positive integers $n$, \begin{equation}
-n\tau=m+\beta+a_{n},\quad|a_{n}|\le\frac{1}{n^{\rho}}\label{eq:3.24}\end{equation}
 and some $\lambda\in\mathbb{S}(\theta)$ with\begin{equation}
n\theta=m_{1}+\lambda,\quad\lambda=\left\{ n\theta\right\} .\label{eq:3.25}\end{equation}
 For each such $\lambda,\beta,n,m$, we have\begin{eqnarray}
S_{n}(x_{n}(z,s);q) & = & \frac{(-z)^{n}q^{n^{2}(1-s)+\left\lfloor m/2\right\rfloor \left(\tau n+\left\lfloor m/2\right\rfloor \right)}}{(q;q)_{\infty}^{2}(-ze^{-2n\theta\pi i})^{\left\lfloor m/2\right\rfloor }}\nonumber \\
 & \times & \left\{ \Theta\left(-zq^{\chi(m)+\beta}e^{-2\lambda\pi i}\mid q\right)+r_{sw}(n|6)\right\} \label{eq:3.26}\end{eqnarray}
 with\begin{eqnarray}
|r_{sw}(n|6)| & \le & \frac{54(-q^{3};q)_{\infty}^{2}\Theta(|z|\mid\sqrt{q})}{(1-q)^{2}}\nonumber \\
 & \times & \left\{ |z|^{\nu_{n}}q^{\nu_{n}^{2}}+\frac{q^{\nu_{n}^{2}/2}}{|z|^{\nu_{n}}}+\frac{\log^{2}n}{n^{\rho}}\right\} \label{eq:3.27}\end{eqnarray}
 and\begin{equation}
\nu_{n}=\left\lfloor \frac{q^{4}\log^{2}n}{1+\log q^{-1}}\right\rfloor \label{eq:3.28}\end{equation}
 for $n$ sufficiently large. 
\item When $-2<\tau<0$ , both $\tau$ and $\theta$ are irrational. Assume
that there exist two real numbers $\beta_{1},\beta_{2}\in[0,1)$ with
$\omega(-\tau,\theta|\beta_{1},\beta_{2})>1$. Then for any \begin{equation}
0<\rho<\omega(-\tau,\theta|\beta_{1},\beta_{2})-1\label{eq:3.29}\end{equation}
 there are infinitely many positive integers $n$ such that\begin{equation}
-\tau n=m+\beta_{1}+a_{n},\quad|a_{n}|<\frac{1}{n^{\rho}}\label{eq:3.30}\end{equation}
 and\begin{equation}
n\theta=m_{1}+\beta_{2}+b_{n},\quad|b_{n}|<\frac{1}{n^{\rho}}.\label{eq:3.31}\end{equation}
 For each such $\beta_{1},\beta_{2},n,m$, we have\begin{eqnarray}
S_{n}(x_{n}(z,s);q) & = & \frac{(-z)^{n}q^{n^{2}(1-s)+\left\lfloor m/2\right\rfloor \left(\tau n+\left\lfloor m/2\right\rfloor \right)}}{(q;q)_{\infty}^{2}(-ze^{-2n\theta\pi i})^{\left\lfloor m/2\right\rfloor }}\nonumber \\
 & \times & \left\{ \Theta\left(-zq^{\chi(m)+\beta_{1}}e^{-2\beta_{2}\pi i}\mid q\right)+r_{sw}(n|7)\right\} \label{eq:3.32}\end{eqnarray}
 with\begin{eqnarray}
|r_{sw}(n|7)| & \le & \frac{54(-q^{3};q)_{\infty}^{2}\Theta(|z|\mid\sqrt{q})}{(1-q)^{2}}\nonumber \\
 & \times & \left\{ |z|^{\nu_{n}}q^{\nu_{n}^{2}}+\frac{q^{\nu_{n}^{2}/2}}{|z|^{\nu_{n}}}+\frac{\log^{2}n}{n^{\rho}}\right\} \label{eq:3.33}\end{eqnarray}
 and\begin{equation}
\nu_{n}=\left\lfloor \frac{q^{4}\log^{2}n}{1+\log q^{-1}}\right\rfloor \label{eq:3.34}\end{equation}
 for $n$ sufficiently large. 
\end{enumerate}
\end{thm}
\begin{rem}
In the cases when only $\tau$ or $\theta$ is irrational, say $\tau$,
if it is an algebraic number of degree $l$, then we could take the
corresponding $\beta=0$, and the order of the error term is no better
than $\mathcal{O}(n^{1-l})$, when this number is a Liouville number,
however, the order of the error term is better than any $\mathcal{O}(n^{-k})$.
It is clear that the scaling \eqref{eq:2.12} with $\tau\le-2$ won't
give us anything interesting. 
\end{rem}

\section{Stieltjes-Wigert Polynomials\label{sec:Stieltjes-Wigert}}

In this section, we have used the following inequalities \begin{equation}
0<(a;q)_{n}\le1,\quad(-b,q)_{n}\ge1\label{eq:4.1}\end{equation}
 for $0\le a<1$, $b\ge0$ and for $n\in\mathbb{N}$ or $n=\infty$.
In the proofs \ref{sub:4.1}--\ref{sub:4.3} we will use the inequalities\begin{equation}
0<\frac{(q;q)_{\infty}}{(q;q)_{n-k}}\le1\label{eq:4.2}\end{equation}
 for $0\le k\le n$ and \begin{equation}
\left|\frac{(q;q)_{\infty}}{(q;q)_{n-k}}-1\right|\le\frac{(-q^{3};q)_{\infty}q^{n/2}}{(1-q)}\label{eq:4.3}\end{equation}
 for $0\le k\le\left\lfloor \frac{n}{2}\right\rfloor $. This can
be seen from \begin{equation}
\frac{(q;q)_{\infty}}{(q;q)_{n-k}}=R_{1}(q;n-k)+1\label{eq:4.4}\end{equation}
 and Lemma \ref{lem:1}.

In the proofs \ref{sub:4.4}--\ref{sub:4.7}, we have \begin{equation}
-2<\tau<0.\label{eq:4.5}\end{equation}
 and\begin{equation}
2\ll\nu_{n}\le n\min\left\{ \frac{2+\tau}{8},\frac{-\tau}{8}\right\} <\frac{n}{4}\label{eq:4.6}\end{equation}
 for $n$ sufficiently large. Assume that\begin{equation}
-\tau n=m+c_{n},\quad0\le c_{n}<1,\quad0<m=\left\lfloor -\tau n\right\rfloor ,\label{eq:4.7}\end{equation}
 \begin{equation}
n\theta=m_{1}+d_{n},\quad0\le d_{n}<1,\quad m_{1}=\left\lfloor n\theta\right\rfloor .\label{eq:4.8}\end{equation}
 Then\begin{eqnarray}
\frac{S_{n}(x_{n}(z,s);q)}{(-z)^{n}q^{n^{2}(1-s)}} & = & \sum_{k=0}^{n}\frac{q^{k^{2}}e^{2nk\theta\pi i}}{(q;q)_{k}(q;q)_{n-k}}\left(-\frac{q^{\tau n}}{z}\right)^{k}\nonumber \\
 & = & \sum_{k=0}^{\left\lfloor m/2\right\rfloor }\frac{q^{k^{2}}e^{2nk\theta\pi i}}{(q;q)_{k}(q;q)_{n-k}}\left(-\frac{q^{\tau n}}{z}\right)^{k}\nonumber \\
 & + & \sum_{k=\left\lfloor m/2\right\rfloor +1}^{n}\frac{q^{k^{2}}e^{2nk\theta\pi i}}{(q;q)_{k}(q;q)_{n-k}}\left(-\frac{q^{\tau n}}{z}\right)^{k}\nonumber \\
 & = & s_{1}+s_{2}.\label{eq:4.9}\end{eqnarray}
 We reverse the summation order in $s_{1}$ to obtain\begin{equation}
\frac{s_{1}(q;q)_{\infty}^{2}(-ze^{-2n\theta\pi i})^{\left\lfloor m/2\right\rfloor }}{q^{\left\lfloor m/2\right\rfloor (\tau n+\left\lfloor m/2\right\rfloor )}}=\sum_{k=0}^{\left\lfloor m/2\right\rfloor }q^{k^{2}}\left(-zq^{\chi(m)+c_{n}}e^{-2\pi id_{n}}\right)^{k}e(k,n)\label{eq:4.10}\end{equation}
 with\begin{equation}
e(k,n)=\frac{(q;q)_{\infty}^{2}}{(q;q)_{\left\lfloor m/2\right\rfloor -k}(q;q)_{n-\left\lfloor m/2\right\rfloor +k}}.\label{eq:4.11}\end{equation}
 Since $0<q<1$, it is clear that\begin{equation}
|e(k,n)|\le1\label{eq:4.12}\end{equation}
 for $0\le k\le\left\lfloor \frac{m}{2}\right\rfloor $. Observe that\begin{eqnarray}
e(k,n) & = & \left\{ R_{1}(q;\left\lfloor \frac{m}{2}\right\rfloor -k)+1\right\} \left\{ R_{1}(q;n-\left\lfloor \frac{m}{2}\right\rfloor +k)+1\right\} .\label{eq:7.13}\end{eqnarray}
 By Lemma \ref{lem:1}, we get\begin{equation}
|e(k,n)-1|\le\frac{3(-q^{3};q)_{\infty}^{2}q^{\nu_{n}+2}}{(1-q)^{2}}\label{eq:4.14}\end{equation}
 for $0\le k\le\nu_{n}-1$.

In sum $s_{2}$ we shift summation from $k$ to $k+\left\lfloor \frac{m}{2}\right\rfloor $
to get \begin{equation}
\frac{s_{2}(q;q)_{\infty}^{2}(-ze^{-2n\theta\pi i})^{\left\lfloor m/2\right\rfloor }}{q^{\left\lfloor m/2\right\rfloor (\tau n+\left\lfloor m/2\right\rfloor )}}=\sum_{k=1}^{n-\left\lfloor m/2\right\rfloor }q^{k^{2}}\left(-z^{-1}q^{-\chi(m)-c_{n}}e^{2\pi id_{n}}\right)^{k}f(k,n)\label{eq:4.15}\end{equation}
 with\begin{equation}
f(k,n)=\frac{(q;q)_{\infty}^{2}}{(q;q)_{\left\lfloor m/2\right\rfloor +k}(q;q)_{n-\left\lfloor m/2\right\rfloor -k}}\label{eq:4.16}\end{equation}
 for $1\le k\le n-\left\lfloor \frac{m}{2}\right\rfloor $. Hence\begin{equation}
|f(k,n)|\le1\label{eq:4.17}\end{equation}
 for $1\le k\le n-\left\lfloor \frac{m}{2}\right\rfloor $. Apply
Lemma \ref{lem:1} to\begin{eqnarray}
f(k,n) & = & \left\{ R_{1}(q;\left\lfloor \frac{m}{2}\right\rfloor +k)+1\right\} \left\{ R_{1}(q;n-\left\lfloor \frac{m}{2}\right\rfloor -k)+1\right\} ,\label{eq:4.18}\end{eqnarray}
 we obtain \begin{equation}
|f(k,n)-1|\le\frac{3(-q^{3};q)_{\infty}^{2}q^{\nu_{n}+2}}{(1-q)^{2}}\label{eq:4.19}\end{equation}
 for $1\le k\le\nu_{n}-1$.

\subsection{Proof for the case $\tau>0$ \label{sub:4.1}}

\begin{proof}
From \eqref{eq:2.13} one has\begin{equation}
\frac{S_{n}(x_{n}(z,s);q)(q;q)_{\infty}}{(-z)^{n}q^{n^{2}(1-s)}}=1+r_{sw}(n|1),\label{eq:4.20}\end{equation}
 and\begin{equation}
r_{sw}(n|1)=\sum_{k=1}^{n}\frac{q^{k^{2}}e^{2nk\theta\pi i}}{(q;q)_{k}}\left(-\frac{q^{\tau n}}{z}\right)^{k}\frac{(q;q)_{\infty}}{(q;q)_{n-k}}.\label{eq:4.21}\end{equation}
 Thus\begin{eqnarray}
|r_{sw}(n|1)| & \le & \sum_{k=1}^{n}\frac{q^{k^{2}}}{(q;q)_{k}}\frac{(q;q)_{\infty}}{(q;q)_{n-k}}\left(\frac{q^{\tau n}}{\left|z\right|}\right)^{k}\nonumber \\
 & \le & \sum_{k=1}^{n}\frac{q^{k^{2}}}{(q;q)_{k}}\left(\frac{q^{\tau n}}{\left|z\right|}\right)^{k}\nonumber \\
 & \le & \sum_{k=1}^{\infty}\frac{q^{k^{2}}}{(q;q)_{k}}\left(\frac{q^{\tau n}}{\left|z\right|}\right)^{k}\nonumber \\
 & \le & \frac{qB_{q}(q^{2}|z|^{-1})}{(1-q)|z|}q^{\tau n}.\label{eq:4.22}\end{eqnarray}

\end{proof}

\subsection{Proof for the case $\tau=0$ and $\theta$ rational\label{sub:4.2}}

\begin{proof}
We have\begin{eqnarray}
\frac{S_{n}(x_{n}(z,s);q)(q;q)_{\infty}}{(-z)^{n}q^{n^{2}(1-s)}} & = & \sum_{k=0}^{n}\frac{q^{k^{2}}}{(q;q)_{k}(q;q)_{n-k}}\left(-\frac{e^{2\lambda\pi i}}{z}\right)^{k}\nonumber \\
 & = & \sum_{k=0}^{\infty}\frac{q^{k^{2}}}{(q;q)_{k}}\left(-\frac{e^{2\lambda\pi i}}{z}\right)^{k}\nonumber \\
 & - & \sum_{k=\left\lfloor n/2\right\rfloor }^{\infty}\frac{q^{k^{2}}}{(q;q)_{k}}\left(-\frac{e^{2\lambda\pi i}}{z}\right)^{k}\nonumber \\
 & + & \sum_{k=0}^{\left\lfloor n/2\right\rfloor -1}\frac{q^{k^{2}}}{(q;q)_{k}}\left(-\frac{e^{2\lambda\pi i}}{z}\right)^{k}\left\{ \frac{(q;q)_{\infty}}{(q;q)_{n-k}}-1\right\} \nonumber \\
 & + & \sum_{k=\left\lfloor n/2\right\rfloor }^{n}\frac{q^{k^{2}}}{(q;q)_{k}}\left(-\frac{e^{2\lambda\pi i}}{z}\right)^{k}\frac{(q;q)_{\infty}}{(q;q)_{n-k}}\nonumber \\
 & = & A_{q}\left(\frac{e^{2\lambda\pi i}}{z}\right)+s_{1}+s_{2}+s_{3}.\label{eq:4.23}\end{eqnarray}
 Then\begin{eqnarray}
|s_{1}+s_{3}| & \le & 2\sum_{k=\left\lfloor n/2\right\rfloor }^{\infty}\frac{q^{k^{2}}}{(q;q)_{k}}\left(\frac{1}{\left|z\right|}\right)^{k}\nonumber \\
 & = & 2\frac{q^{n^{2}/4}B_{q}(|z|^{-1})}{(q;q)_{\infty}|z|^{\left\lfloor n/2\right\rfloor }},\label{eq:4.24}\end{eqnarray}
 and \begin{equation}
|s_{2}|\le\frac{(-q^{3};q)_{\infty}B_{q}(|z|^{-1})q^{n/2}}{(1-q)}.\label{eq:4.25}\end{equation}
 Therefore \begin{equation}
\frac{S_{n}(x_{n}(z,s);q)(q;q)_{\infty}}{(-z)^{n}q^{n^{2}(1-s)}}=A_{q}\left(\frac{e^{2\lambda\pi i}}{z}\right)+r_{sw}(n|2)\label{eq:4.26}\end{equation}
 with\begin{equation}
r_{sw}(n|2)=s_{1}+s_{2}+s_{3}\label{eq:4.27}\end{equation}
 and\begin{equation}
|r_{sw}(n|2)|\le\frac{2(-q^{3};q)_{\infty}B_{q}(|z|^{-1})}{(q;q)_{\infty}}\left\{ q^{n/2}+\frac{q^{n^{2}/4}}{|z|^{\left\lfloor n/2\right\rfloor }}\right\} .\label{eq:4.28}\end{equation}

\end{proof}

\subsection{Proof for the case $\tau=0$ and $\theta$ irrational\label{sub:4.3}}

\begin{proof}
Assume that \begin{equation}
\nu_{n}=\left\lfloor \frac{q^{4}\log^{2}n}{1+\log q^{-1}}\right\rfloor ,\label{eq:4.29}\end{equation}
 we have\begin{eqnarray}
\frac{S_{n}(x_{n}(z,s);q)}{(-z)^{n}q^{n^{2}(1-s)}} & = & \sum_{k=0}^{n}\frac{q^{k^{2}}e^{2k\gamma_{n}\pi i}}{(q;q)_{k}(q;q)_{n-k}}\left(-\frac{e^{2\pi i\beta}}{z}\right)^{k}\nonumber \\
 & = & \sum_{k=0}^{\infty}\frac{q^{k^{2}}}{(q;q)_{k}}\left(-\frac{e^{2\pi i\beta}}{z}\right)^{k}\nonumber \\
 & - & \sum_{k=\nu_{n}}^{\infty}\frac{q^{k^{2}}}{(q;q)_{k}}\left(-\frac{e^{2\pi i\beta}}{z}\right)^{k}\nonumber \\
 & + & \sum_{k=0}^{\nu_{n}-1}\frac{q^{k^{2}}}{(q;q)_{k}}\left(-\frac{e^{2\pi i\beta}}{z}\right)^{k}\left\{ \frac{(q;q)_{\infty}}{(q;q)_{n-k}}-1\right\} \nonumber \\
 & + & \sum_{k=0}^{\nu_{n}-1}\frac{q^{k^{2}}}{(q;q)_{k}}\left(-\frac{e^{2\pi i\beta}}{z}\right)^{k}\frac{(q;q)_{\infty}}{(q;q)_{n-k}}\left\{ e^{2\pi ik\gamma_{n}}-1\right\} \nonumber \\
 & + & \sum_{k=\nu_{n}}^{n}\frac{q^{k^{2}}}{(q;q)_{k}}\left(-\frac{e^{2\pi i\beta}}{z}\right)^{k}\frac{(q;q)_{\infty}}{(q;q)_{n-k}}e^{2\pi ik\gamma_{n}}\nonumber \\
 & = & A_{q}\left(\frac{e^{2\pi i\beta}}{z}\right)+s_{1}+s_{2}+s_{3}+s_{4}.\label{eq:4.30}\end{eqnarray}
 Then,\begin{eqnarray}
|s_{1}+s_{4}| & \le & 2\sum_{k=\nu_{n}}^{\infty}\frac{q^{k^{2}}|z|^{-k}}{(q;q)_{k}}\nonumber \\
 & \le & \frac{2B_{q}(|z|^{-1})q^{\nu_{n}^{2}}}{(q;q)_{\infty}|z|^{\nu_{n}}}\label{eq:4.31}\end{eqnarray}
 for $n$ sufficiently large.

Clearly, we also have\begin{equation}
\nu_{n}\ll n^{\min(1,\rho)}/8\label{eq:4.32}\end{equation}
 and\begin{equation}
q^{n/2}\ll\frac{\nu_{n}}{n^{\rho}}\label{eq:4.33}\end{equation}
 for sufficiently large $n$. Hence, \begin{equation}
|s_{2}|\le\frac{(-q^{3};q)_{\infty}B_{q}(|z|^{-1})}{(1-q)}\frac{\log^{2}n}{n^{\rho}}\label{eq:4.34}\end{equation}
 and\begin{equation}
|e^{2\pi ik\gamma_{n}}-1|\le\frac{2\pi\nu_{n}}{n^{\rho}}e^{2\pi\nu_{n}/n^{\rho}}\le\frac{24\log^{2}n}{n^{\rho}},\label{eq:4.35}\end{equation}
 and \begin{equation}
|s_{3}|\le\frac{24B_{q}(|z|^{-1})\log^{2}n}{n^{\rho}}\label{eq:4.35}\end{equation}
 for $n$ sufficiently large.

Let \begin{equation}
e_{sw}(n|3)=s_{1}+s_{2}+s_{3}+s_{4},\label{eq:4.36}\end{equation}
 then,\begin{equation}
\frac{S_{n}(x_{n}(z,s);q)}{(-z)^{n}q^{n^{2}(1-s)}}=A_{q}\left(\frac{e^{2\pi i\beta}}{z}\right)+e_{sw}(n|3)\label{eq:4.37}\end{equation}
 with\begin{equation}
|e_{sw}(n|3)|\le\frac{24(-q^{3};q)_{\infty}B_{q}(|z|^{-1})}{(q;q)_{\infty}}\left\{ \frac{\log^{2}n}{n^{\rho}}+\frac{q^{\nu_{n}^{2}}}{|z|^{\nu_{n}}}\right\} \label{eq:4.38}\end{equation}
 for $n$ sufficiently large.
\end{proof}

\subsection{Proof for the case $-2<\tau<0$ rational and $\theta$ rational\label{sub:4.4}}

\begin{proof}
Assume that\begin{equation}
\nu_{n}=\min\left\{ \left\lfloor \frac{(2+\tau)n}{8}\right\rfloor ,\left\lfloor \frac{-\tau n}{8}\right\rfloor \right\} .\label{eq:4.39}\end{equation}
 From \eqref{eq:4.10} we get\begin{eqnarray}
\frac{s_{1}(q;q)_{\infty}^{2}(-ze^{-2n\theta\pi i})^{\left\lfloor m/2\right\rfloor }}{q^{\left\lfloor m/2\right\rfloor (\tau n+\left\lfloor m/2\right\rfloor )}} & = & \sum_{k=0}^{\left\lfloor m/2\right\rfloor }q^{k^{2}}\left(-zq^{\chi(m)+\lambda}e^{-2\pi i\lambda_{1}}\right)^{k}e(k,n)\nonumber \\
 & = & \sum_{k=0}^{\infty}q^{k^{2}}\left(-zq^{\chi(m)+\lambda}e^{-2\pi i\lambda_{1}}\right)^{k}\nonumber \\
 & - & \sum_{k=\nu_{n}}^{\infty}q^{k^{2}}\left(-zq^{\chi(m)+\lambda}e^{-2\pi i\lambda_{1}}\right)^{k}\nonumber \\
 & + & \sum_{k=0}^{\nu_{n}-1}q^{k^{2}}\left(-zq^{\chi(m)+\lambda}e^{-2\pi i\lambda_{1}}\right)^{k}\left(e(k,n)-1\right)\nonumber \\
 & + & \sum_{k=\nu_{n}}^{\left\lfloor m/2\right\rfloor }q^{k^{2}}\left(-zq^{\chi(m)+\lambda}e^{-2\pi i\lambda_{1}}\right)^{k}e(k,n)\nonumber \\
 & = & \sum_{k=0}^{\infty}q^{k^{2}}\left(-zq^{\chi(m)+\lambda}e^{-2\pi i\lambda_{1}}\right)^{k}+s_{11}+s_{12}+s_{13}.\label{eq:4.40}\end{eqnarray}
 Then,\begin{eqnarray}
|s_{11}+s_{12}| & \le & 2\sum_{k=\nu_{n}}^{\infty}q^{k^{2}}\left(\left|z\right|q^{\chi(m)+\lambda}\right)^{k}\nonumber \\
 & \le2 & |z|^{\nu_{n}}q^{\nu_{n}^{2}}\sum_{k=0}^{\infty}q^{k^{2}}\left(\left|z\right|q^{\chi(m)+\lambda}\right)^{k}\nonumber \\
 & \le & 2\Theta\left(\left|z\right|\mid\sqrt{q}\right)q^{\nu_{n}^{2}}|z|^{\nu_{n}}\label{eq:4.41}\end{eqnarray}
 and\begin{eqnarray}
|s_{12}| & \le & \frac{3(-q^{3};q)_{\infty}^{2}q^{\nu_{n}}}{(1-q)^{2}}\sum_{k=0}^{\infty}q^{k^{2}}(|z|q^{\chi(m)+\lambda})^{k}\nonumber \\
 & \le & \frac{3(-q^{3};q)_{\infty}^{2}\Theta(|z|\mid\sqrt{q})q^{\nu_{n}}}{(1-q)^{2}}.\label{eq:4.42}\end{eqnarray}
 Let \begin{equation}
r_{1}(n)=s_{11}+s_{12}+s_{13},\label{eq:4.43}\end{equation}
 then,\begin{equation}
\frac{s_{1}(q;q)_{\infty}^{2}(-ze^{-2n\theta\pi i})^{\left\lfloor m/2\right\rfloor }}{q^{\left\lfloor m/2\right\rfloor (\tau n+\left\lfloor m/2\right\rfloor )}}=\sum_{k=0}^{\infty}q^{k^{2}}\left(-zq^{\chi(m)+\lambda}e^{-2\pi i\lambda_{1}}\right)^{k}+r_{1}(n)\label{eq:4.44}\end{equation}
 \begin{equation}
|r_{1}(n)|\le\frac{3(-q^{3};q)_{\infty}^{2}\Theta(|z|\mid\sqrt{q})}{(1-q)^{2}}\left\{ q^{\nu_{n}^{2}}|z|^{\nu_{n}}+q^{\nu_{n}}\right\} .\label{eq:4.45}\end{equation}
 From \eqref{eq:4.15} we obtain,\begin{eqnarray}
\frac{s_{2}(q;q)_{\infty}^{2}(-ze^{-2n\theta\pi i})^{\left\lfloor m/2\right\rfloor }}{q^{\left\lfloor m/2\right\rfloor (\tau n+\left\lfloor m/2\right\rfloor )}} & = & \sum_{k=1}^{n-\left\lfloor m/2\right\rfloor }q^{k^{2}}\left(-z^{-1}q^{-\chi(m)-\lambda}e^{2\pi i\lambda_{1}}\right)^{k}f(k,n)\nonumber \\
 & = & \sum_{k=1}^{\infty}q^{k^{2}}\left(-z^{-1}q^{-\chi(m)-\lambda}e^{2\pi i\lambda_{1}}\right)^{k}\nonumber \\
 & - & \sum_{k=\nu_{n}}^{\infty}q^{k^{2}}\left(-z^{-1}q^{-\chi(m)-\lambda}e^{2\pi i\lambda_{1}}\right)^{k}\nonumber \\
 & + & \sum_{k=1}^{\nu_{n}-1}q^{k^{2}}\left(-z^{-1}q^{-\chi(m)-\lambda}e^{2\pi i\lambda_{1}}\right)^{k}\left(f(k,n)-1\right)\nonumber \\
 & + & \sum_{k=\nu_{n}}^{n-\left\lfloor m/2\right\rfloor }q^{k^{2}}\left(-z^{-1}q^{-\chi(m)-\lambda}e^{2\pi i\lambda_{1}}\right)^{k}f(k,n)\nonumber \\
 & = & \sum_{k=-1}^{-\infty}q^{k^{2}}\left(-zq^{\chi(m)+\lambda}e^{-2\pi i\lambda_{1}}\right)^{k}+s_{21}+s_{22}+s_{23},\label{eq:4.46}\end{eqnarray}
 then, \begin{eqnarray}
|s_{21}+s_{23}| & \le & 2\sum_{k=\nu_{n}}^{\infty}q^{k^{2}}\left(\left|z\right|^{-1}q^{-\chi(m)-\lambda}\right)^{k}\nonumber \\
 & \le & 2\frac{q^{\nu_{n}^{2}-2\nu_{n}}\sum_{k=0}^{\infty}q^{k^{2}}|z|^{-k}}{{|z|}^{\nu_{n}}}\nonumber \\
 & \le & 2\frac{q^{\nu_{n}^{2}/2}\sum_{k=0}^{\infty}q^{k^{2}/2}|z|^{-k}}{{|z|}^{\nu_{n}}}\nonumber \\
 & \le & 2\frac{q^{\nu_{n}^{2}/2}\Theta(|z|\mid\sqrt{q})}{{|z|}^{\nu_{n}}}\label{eq:4.47}\end{eqnarray}
 for $n$ sufficiently large. For the sum $s_{22}$, we have\begin{eqnarray}
|s_{22}| & \le & \frac{3(-q^{3};q)_{\infty}^{2}q^{\nu_{n}+2}}{(1-q)^{2}}\sum_{k=1}^{\infty}q^{k^{2}/2+k^{2}/2-2k}|z|^{-k}\nonumber \\
 & \le & \frac{3(-q^{3};q)_{\infty}^{2}q^{\nu_{n}}}{(1-q)^{2}}\sum_{k=1}^{\infty}q^{k^{2}/2}|z|^{-k}\nonumber \\
 & \le & \frac{3(-q^{3};q)_{\infty}^{2}\Theta(|z|\mid\sqrt{q})q^{\nu_{n}}}{(1-q)^{2}}.\label{eq:4.48}\end{eqnarray}
 Hence \begin{equation}
\frac{s_{2}(q;q)_{\infty}^{2}(-ze^{-2n\theta\pi i})^{\left\lfloor m/2\right\rfloor }}{q^{\left\lfloor m/2\right\rfloor (\tau n+\left\lfloor m/2\right\rfloor )}}=\sum_{k=-1}^{-\infty}q^{k^{2}}\left(-zq^{\chi(m)+\lambda}e^{-2\pi i\lambda_{1}}\right)^{k}+r_{2}(n)\label{eq:4.49}\end{equation}
 with\begin{equation}
r_{2}(n)=s_{21}+s_{22}+s_{23}\label{eq:4.50}\end{equation}
 and\begin{equation}
|r_{2}(n)|\le\frac{3(-q^{3};q)_{\infty}^{2}\Theta(|z|\mid\sqrt{q})}{(1-q)^{2}}\left\{ q^{\nu_{n}}+\frac{q^{\nu_{n}^{2}/2}}{{|z|}^{\nu_{n}}}\right\} .\label{eq:4.51}\end{equation}
 Thus\begin{equation}
\frac{(-ze^{-2n\theta\pi i})^{\left\lfloor m/2\right\rfloor }S_{n}(x_{n}(z,s);q)(q;q)_{\infty}^{2}}{(-z)^{n}q^{n^{2}(1-s)+\left\lfloor m/2\right\rfloor \left(\tau n+\left\lfloor m/2\right\rfloor \right)}}=\Theta\left(-zq^{\chi(m)+\lambda}e^{-2\pi i\lambda_{1}}\mid q\right)+r_{sw}(n|4)\label{eq:4.52}\end{equation}
 with\begin{equation}
|r_{sw}(n|4)|\le\frac{6(-q^{3};q)_{\infty}^{2}\Theta(|z|\mid\sqrt{q})}{(1-q)^{2}}\left\{ q^{\nu_{n}}+q^{\nu_{n}^{2}}|z|^{\nu_{n}}+\frac{q^{\nu_{n}^{2}/2}}{{|z|}^{\nu_{n}}}\right\} \label{eq:4.53}\end{equation}
 for $n$ sufficiently large.
\end{proof}

\subsection{Proof for the case $-2<\tau<0$ rational and $\theta$ irrational\label{sub:4.5}}

\begin{proof}
The existence of $\lambda$ with infinitely many positive integers
$n$ satisfying \eqref{eq:3.18} and \eqref{eq:3.19} is guaranteed
by the fact that $\mathbb{S}(-\tau)$ is a finite set.

It is clear that for sufficiently large $n$, we have\begin{equation}
{2\ll\nu}_{n}=\left\lfloor \frac{q^{4}\log^{2}n}{1+\log q^{-1}}\right\rfloor \ll\frac{n^{\rho}}{8},\quad q^{\nu_{n}}\ll\frac{\nu_{n}}{n^{\rho}}.\label{eq:4.54}\end{equation}
 From \eqref{eq:4.10} ,\begin{eqnarray}
\frac{s_{1}(q;q)_{\infty}^{2}(-ze^{-2n\theta\pi i})^{\left\lfloor m/2\right\rfloor }}{q^{\left\lfloor m/2\right\rfloor (\tau n+\left\lfloor m/2\right\rfloor )}} & = & \sum_{k=0}^{\left\lfloor m/2\right\rfloor }q^{k^{2}}\left(-zq^{\chi(m)+\lambda}e^{-2\pi i\beta-2\pi ib_{n}}\right)^{k}e(k,n)\nonumber \\
 & = & \sum_{k=0}^{\infty}q^{k^{2}}\left(-zq^{\chi(m)+\lambda}e^{-2\pi i\beta}\right)^{k}\nonumber \\
 & - & \sum_{k=\nu_{n}}^{\infty}q^{k^{2}}\left(-zq^{\chi(m)+\lambda}e^{-2\pi i\beta}\right)^{k}\nonumber \\
 & + & \sum_{k=0}^{\nu_{n}-1}q^{k^{2}}\left(-zq^{\chi(m)+\lambda}e^{-2\pi i\beta}\right)^{k}\left\{ e^{-2k\pi ib_{n}}-1\right\} \nonumber \\
 & + & \sum_{k=0}^{\nu_{n}-1}q^{k^{2}}\left(-zq^{\chi(m)+\lambda}e^{-2\pi i\beta}\right)^{k}e^{-2k\pi ib_{n}}\left\{ e(k,n)-1\right\} \nonumber \\
 & + & \sum_{k=\nu_{n}}^{\left\lfloor m/2\right\rfloor }q^{k^{2}}\left(-zq^{\chi(m)+\lambda}e^{-2\pi i\beta}\right)^{k}e^{-2k\pi ib_{n}}e(k,n)\nonumber \\
 & = & \sum_{k=0}^{\infty}q^{k^{2}}\left(-zq^{\chi(m)+\lambda}e^{-2\pi i\beta}\right)^{k}+s_{11}+s_{12}+s_{13}+s_{14},\label{eq:4.55}\end{eqnarray}
 then\begin{eqnarray}
|s_{11}+s_{14}| & \le & 2\sum_{k=\nu_{n}}^{\infty}q^{k^{2}}{|z|}^{k}\nonumber \\
 & \le & 2|z|^{\nu_{n}}q^{\nu_{n}^{2}}\sum_{k=0}^{\infty}q^{k^{2}}{|z|}^{k}\nonumber \\
 & \le & 2\Theta(|z|\mid\sqrt{q})|z|^{\nu_{n}}q^{\nu_{n}^{2}}.\label{eq:4.56}\end{eqnarray}
 For sufficiently large $n$, we have\begin{eqnarray}
|s_{12}| & \le & \frac{24\nu_{n}}{n^{\rho}}\sum_{k=0}^{\infty}q^{k^{2}}{|z|}^{k}\nonumber \\
 & \le & \frac{24\Theta(|z|\mid\sqrt{q})\log^{2}n}{n^{\rho}}\label{eq:4.57}\end{eqnarray}
 and\begin{eqnarray}
|s_{13}| & \le & \frac{3(-q^{3};q)_{\infty}^{2}q^{\nu_{n}+2}}{(1-q)^{2}}\sum_{k=0}^{\infty}q^{k^{2}}{|z|}^{k}\nonumber \\
 & \le & \frac{3(-q^{3};q)_{\infty}^{2}\Theta(|z|\mid\sqrt{q})}{(1-q)^{2}}\frac{\log^{2}n}{n^{\rho}}.\label{eq:4.58}\end{eqnarray}
 Hence, \begin{equation}
\frac{s_{1}(q;q)_{\infty}^{2}(-ze^{-2n\theta\pi i})^{\left\lfloor m/2\right\rfloor }}{q^{\left\lfloor m/2\right\rfloor (\tau n+\left\lfloor m/2\right\rfloor )}}=\sum_{k=0}^{\infty}q^{k^{2}}\left(-zq^{\chi(m)+\lambda}e^{-2\pi i\beta}\right)^{k}+r_{1}(n)\label{eq:4.59}\end{equation}
 with\begin{equation}
r_{1}(n)=s_{11}+s_{12}+s_{13}+s_{14}\label{eq:4.60}\end{equation}
 and\begin{equation}
|r_{1}(n)|\le\frac{27(-q^{3};q)_{\infty}^{2}\Theta(|z|\mid\sqrt{q})}{(1-q)^{2}}\left\{ |z|^{\nu_{n}}q^{\nu_{n}^{2}}+\frac{\log^{2}n}{n^{\rho}}\right\} \label{eq:4.61}\end{equation}
 for $n$ sufficiently large.

From \eqref{eq:4.15} we obtain

\begin{eqnarray}
\frac{s_{2}(q;q)_{\infty}^{2}(-ze^{-2n\theta\pi i})^{\left\lfloor m/2\right\rfloor }}{q^{\left\lfloor m/2\right\rfloor (\tau n+\left\lfloor m/2\right\rfloor )}} & = & \sum_{k=1}^{n-\left\lfloor m/2\right\rfloor }q^{k^{2}}\left(-z^{-1}q^{-\chi(m)-\lambda}e^{2\beta\pi i}\right)^{k}e^{2k\pi ib_{n}}f(k,n)\nonumber \\
 & = & \sum_{k=1}^{\infty}q^{k^{2}}\left(-z^{-1}q^{-\chi(m)-\lambda}e^{2\beta\pi i}\right)^{k}\nonumber \\
 & - & \sum_{k=\nu_{n}}^{\infty}q^{k^{2}}\left(-z^{-1}q^{-\chi(m)-\lambda}e^{2\beta\pi i}\right)^{k}\nonumber \\
 & + & \sum_{k=1}^{\nu_{n}-1}q^{k^{2}}\left(-z^{-1}q^{-\chi(m)-\lambda}e^{2\beta\pi i}\right)^{k}\left\{ e^{2k\pi ib_{n}}-1\right\} \nonumber \\
 & + & \sum_{k=1}^{\nu_{n}-1}q^{k^{2}}\left(-z^{-1}q^{-\chi(m)-\lambda}e^{2\beta\pi i}\right)^{k}e^{2k\pi ib_{n}}\left\{ f(k,n)-1\right\} \nonumber \\
 & + & \sum_{k=\nu_{n}}^{n-\left\lfloor m/2\right\rfloor }q^{k^{2}}\left(-z^{-1}q^{-\chi(m)-\lambda}e^{2\beta\pi i}\right)^{k}e^{2k\pi ib_{n}}f(k,n)\nonumber \\
 & = & \sum_{k=-1}^{-\infty}q^{k^{2}}\left(-zq^{\chi(m)+\lambda}e^{-2\beta\pi i}\right)^{k}+s_{21}+s_{22}+s_{23}+s_{24},\label{eq:4.62}\end{eqnarray}
 then\begin{eqnarray}
|s_{21}+s_{24}| & \le & 2\sum_{k=\nu_{n}}^{\infty}q^{k^{2}-2k}|z|^{-k}\nonumber \\
 & \le & 2|z|^{-\nu_{n}}q^{\nu_{n}^{2}-2\nu_{n}}\sum_{k=0}^{\infty}q^{k^{2}/2}|z|^{-k}\nonumber \\
 & \le & \frac{2\Theta(|z|\mid\sqrt{q})q^{\nu_{n}^{2}/2}}{|z|^{\nu_{n}}}\label{eq:4.63}\end{eqnarray}
 and \begin{eqnarray}
|s_{22}| & \le & \frac{24\nu_{n}}{n^{\rho}}\sum_{k=0}^{\infty}q^{k^{2}-2k}{|z|}^{-k}\nonumber \\
 & \le & \frac{24\Theta(|z|\mid\sqrt{q})\log^{2}n}{n^{\rho}}\label{eq:4.64}\end{eqnarray}
 and\begin{eqnarray}
|s_{23}| & \le & \frac{3(-q^{3};q)_{\infty}^{2}q^{\nu_{n}+2}}{(1-q)^{2}}\sum_{k=0}^{\infty}q^{k^{2}/2+k^{2}/2-2k}{|z|}^{-k}\nonumber \\
 & \le & \frac{3(-q^{3};q)_{\infty}^{2}q^{\nu_{n}}}{(1-q)^{2}}\sum_{k=0}^{\infty}q^{k^{2}/2}{|z|}^{-k}\nonumber \\
 & \le & \frac{3(-q^{3};q)_{\infty}^{2}\Theta(|z|\mid\sqrt{q})}{(1-q)^{2}}\frac{\log^{2}n}{n^{\rho}}.\label{eq:4.65}\end{eqnarray}
 for $n$ sufficiently large. Therefore\begin{equation}
\frac{s_{2}(q;q)_{\infty}^{2}(-ze^{-2n\theta\pi i})^{\left\lfloor m/2\right\rfloor }}{q^{\left\lfloor m/2\right\rfloor (\tau n+\left\lfloor m/2\right\rfloor )}}=\sum_{k=-1}^{-\infty}q^{k^{2}}\left(-zq^{\chi(m)+\lambda}e^{-2\beta\pi i}\right)^{k}+r_{2}(n)\label{eq:4.66}\end{equation}
 with\begin{equation}
r_{2}(n)=s_{21}+s_{22}+s_{23}+s_{24}\label{eq:4.67}\end{equation}
 and\begin{equation}
|r_{2}(n)|\le\frac{27(-q^{3};q)_{\infty}^{2}\Theta(|z|\mid\sqrt{q})}{(1-q)^{2}}\left\{ \frac{q^{\nu_{n}^{2}/2}}{|z|^{\nu_{n}}}+\frac{\log^{2}n}{n^{\rho}}\right\} \label{eq:4.68}\end{equation}
 for $n$ sufficiently large. Combine \eqref{eq:4.59} and \eqref{eq:4.66}
we get \begin{equation}
\frac{(-ze^{-2n\theta\pi i})^{\left\lfloor m/2\right\rfloor }S_{n}(x_{n}(z,s);q)(q;q)_{\infty}^{2}}{(-z)^{n}q^{n^{2}(1-s)+\left\lfloor m/2\right\rfloor \left(\tau n+\left\lfloor m/2\right\rfloor \right)}}=\Theta\left(-zq^{\chi(m)+\lambda}e^{-2\beta\pi i}\mid q\right)+r_{sw}(n|5)\label{eq:4.69}\end{equation}
 with\begin{equation}
r_{sw}(n|5)=r_{1}(n)+r_{2}(n)\label{eq:4.70}\end{equation}
 and\begin{equation}
|r_{sw}(n|5)|\le\frac{54(-q^{3};q)_{\infty}^{2}\Theta(|z|\mid\sqrt{q})}{(1-q)^{2}}\left\{ |z|^{\nu_{n}}q^{\nu_{n}^{2}}+\frac{q^{\nu_{n}^{2}/2}}{|z|^{\nu_{n}}}+\frac{\log^{2}n}{n^{\rho}}\right\} \label{eq:4.71}\end{equation}
 for $n$ sufficiently large.
\end{proof}

\subsection{Proof for the case $-2<\tau<0$ irrational and $\theta$ rational\label{sub:4.6}}

\begin{proof}
{[}Proof of $-1<\tau<0$ negative irrational and $\theta$ rational]This
proof is very similar to the proof when $\tau$ is rational while
$\theta$ is irrational. We just observe that\begin{align}
 & \frac{s_{1}(q;q)_{\infty}^{2}(-ze^{-2n\theta\pi i})^{\left\lfloor m/2\right\rfloor }}{q^{\left\lfloor m/2\right\rfloor (\tau n+\left\lfloor m/2\right\rfloor )}}\label{eq:4.72}\\
= & \sum_{k=0}^{\left\lfloor m/2\right\rfloor }q^{k^{2}}\left(-zq^{\chi(m)+\beta+a_{n}}e^{-2\pi i\lambda}\right)^{k}e(k,n)\nonumber \\
= & \sum_{k=0}^{\infty}q^{k^{2}}\left(-zq^{\chi(m)+\beta}e^{-2\pi i\lambda}\right)^{k}\nonumber \\
- & \sum_{k=\nu_{n}}^{\infty}q^{k^{2}}\left(-zq^{\chi(m)+\beta}e^{-2\pi i\lambda}\right)^{k}\nonumber \\
+ & \sum_{k=0}^{\nu_{n}-1}q^{k^{2}}\left(-zq^{\chi(m)+\beta}e^{-2\pi i\lambda}\right)^{k}\left\{ q^{ka_{n}}-1\right\} \nonumber \\
+ & \sum_{k=0}^{\nu_{n}-1}q^{k^{2}}\left(-zq^{\chi(m)+\beta}e^{-2\pi i\lambda}\right)^{k}q^{ka_{n}}\left\{ e(k,n)-1\right\} \nonumber \\
+ & \sum_{k=\nu_{n}}^{\left\lfloor m/2\right\rfloor }q^{k^{2}}\left(-zq^{\chi(m)+\beta}e^{-2\pi i\lambda}\right)^{k}q^{ka_{n}}e(k,n),\nonumber \end{align}
 and\begin{align}
 & \frac{s_{2}(q;q)_{\infty}^{2}(-ze^{-2n\theta\pi i})^{\left\lfloor m/2\right\rfloor }}{q^{\left\lfloor m/2\right\rfloor (\tau n+\left\lfloor m/2\right\rfloor )}}\label{eq:4.73}\\
= & \sum_{k=1}^{n-\left\lfloor m/2\right\rfloor }q^{k^{2}}\left(-z^{-1}q^{-\chi(m)-\beta}e^{2\lambda\pi i}\right)^{k}q^{-ka_{n}}f(k,n)\nonumber \\
= & \sum_{k=-1}^{-\infty}q^{k^{2}}\left(-zq^{\chi(m)+\beta}e^{-2\pi i\lambda}\right)^{k}\nonumber \\
- & \sum_{k=\nu_{n}}^{\infty}q^{k^{2}}\left(-z^{-1}q^{-\chi(m)-\beta}e^{2\lambda\pi i}\right)^{k}\nonumber \\
+ & \sum_{k=1}^{\nu_{n}-1}q^{k^{2}}\left(-z^{-1}q^{-\chi(m)-\beta}e^{2\lambda\pi i}\right)^{k}\left\{ q^{-ka_{n}}-1\right\} \nonumber \\
+ & \sum_{k=1}^{\nu_{n}-1}q^{k^{2}}\left(-z^{-1}q^{-\chi(m)-\beta}e^{2\lambda\pi i}\right)^{k}q^{-ka_{n}}\left\{ f(k,n)-1\right\} \nonumber \\
+ & \sum_{k=\nu_{n}}^{n-\left\lfloor m/2\right\rfloor }q^{k^{2}}\left(-z^{-1}q^{-\chi(m)-\beta}e^{2\lambda\pi i}\right)^{k}q^{-ka_{n}}f(k,n).\nonumber \end{align}
 
\end{proof}

\subsection{Proof for the case $-2<\tau<0$ irrational and $\theta$ irrational\label{sub:4.7} }

\begin{proof}
In this case we have\begin{align}
 & \frac{s_{1}(q;q)_{\infty}^{2}(-ze^{-2n\theta\pi i})^{\left\lfloor m/2\right\rfloor }}{q^{\left\lfloor m/2\right\rfloor (\tau n+\left\lfloor m/2\right\rfloor )}}\label{eq:4.74}\\
= & \sum_{k=0}^{\left\lfloor m/2\right\rfloor }q^{k^{2}}\left(-zq^{\chi(m)+\beta_{1}+a_{n}}e^{-2\pi i\beta_{2}-2\pi ib_{n}}\right)^{k}e(k,n)\nonumber \\
= & \sum_{k=0}^{\infty}q^{k^{2}}\left(-zq^{\chi(m)+\beta_{1}}e^{-2\pi i\beta_{2}}\right)^{k}\nonumber \\
- & \sum_{k=\nu_{n}}^{\infty}q^{k^{2}}\left(-zq^{\chi(m)+\beta_{1}}e^{-2\pi i\beta_{2}}\right)^{k}\nonumber \\
+ & \sum_{k=0}^{\nu_{n}-1}q^{k^{2}}\left(-zq^{\chi(m)+\beta_{1}}e^{-2\pi i\beta_{2}}\right)^{k}\left\{ q^{ka_{n}}-1\right\} \nonumber \\
+ & \sum_{k=0}^{\nu_{n}-1}q^{k^{2}}\left(-zq^{\chi(m)+\beta_{1}}e^{-2\pi i\beta_{2}}\right)^{k}\left\{ e^{-2k\pi ib_{n}}-1\right\} q^{ka_{n}}\nonumber \\
+ & \sum_{k=0}^{\nu_{n}-1}q^{k^{2}}\left(-zq^{\chi(m)+\beta_{1}}e^{-2\pi i\beta_{2}}\right)^{k}e^{-2k\pi ib_{n}}q^{ka_{n}}\left\{ e(k,n)-1\right\} \nonumber \\
+ & \sum_{k=\nu_{n}}^{\left\lfloor m/2\right\rfloor }q^{k^{2}}\left(-zq^{\chi(m)+\beta_{1}}e^{-2\pi i\beta_{2}}\right)^{k}e^{-2k\pi ib_{n}}q^{ka_{n}}e(k,n)\nonumber \end{align}

and\begin{align}
 & \frac{s_{2}(q;q)_{\infty}^{2}(-ze^{-2n\theta\pi i})^{\left\lfloor m/2\right\rfloor }}{q^{\left\lfloor m/2\right\rfloor (\tau n+\left\lfloor m/2\right\rfloor )}}\label{eq:4.75}\\
= & \sum_{k=1}^{n-\left\lfloor m/2\right\rfloor }q^{k^{2}}\left(-z^{-1}q^{-\chi(m)-\beta_{1}}e^{2\beta_{2}\pi i}\right)^{k}e^{2k\pi ib_{n}}q^{-ka_{n}}f(k,n)\nonumber \\
= & \sum_{k=-1}^{-\infty}q^{k^{2}}\left(-zq^{\chi(m)+\beta_{1}}e^{-2\beta_{2}\pi i}\right)^{k}\nonumber \\
- & \sum_{k=\nu_{n}}^{\infty}q^{k^{2}}\left(-z^{-1}q^{-\chi(m)-\beta_{1}}e^{2\beta_{2}\pi i}\right)^{k}\nonumber \\
+ & \sum_{k=1}^{\nu_{n}-1}q^{k^{2}}\left(-z^{-1}q^{-\chi(m)-\beta_{1}}e^{2\beta_{2}\pi i}\right)^{k}\left\{ e^{2k\pi ib_{n}}-1\right\} \nonumber \\
+ & \sum_{k=1}^{\nu_{n}-1}q^{k^{2}}\left(-z^{-1}q^{-\chi(m)-\beta_{1}}e^{2\beta_{2}\pi i}\right)^{k}e^{2k\pi ib_{n}}\left\{ q^{-ka_{n}}-1\right\} \nonumber \\
+ & \sum_{k=1}^{\nu_{n}-1}q^{k^{2}}\left(-z^{-1}q^{-\chi(m)-\beta_{1}}e^{2\beta_{2}\pi i}\right)^{k}e^{2k\pi ib_{n}}q^{-ka_{n}}\left\{ f(k,n)-1\right\} \nonumber \\
+ & \sum_{k=\nu_{n}}^{n-\left\lfloor m/2\right\rfloor }q^{k^{2}}\left(-z^{-1}q^{-\chi(m)-\beta_{1}}e^{2\beta_{2}\pi i}\right)^{k}e^{2k\pi ib_{n}}q^{-ka_{n}}f(k,n).\nonumber \end{align}
 To finish the proof, we have to estimate the sums above, but they
are very similar to what we have done in the other cases, we won't
repeat these details here.
\end{proof}

\end{document}